\documentclass[sn-mathphys,Numbered]{sn-jnl}
\usepackage{graphicx}%
\usepackage{multirow}%
\usepackage{amsmath,amssymb,amsfonts}%
\usepackage{amsthm}%
\usepackage{mathrsfs}%
\usepackage[title]{appendix}%
\usepackage{xcolor}%
\usepackage{textcomp}%
\usepackage{manyfoot}%
\usepackage{booktabs}%
\usepackage{algorithm}%
\usepackage{algorithmicx}%
\usepackage{algpseudocode}%
\usepackage{listings}%
\usepackage{cleveref}
\usepackage{natbib}

\usepackage{cuted}

\let\emptyset\varnothing

\newtheorem{theorem}{Theorem}[section]
\newtheorem{corollary}[theorem]{Corollary}
\newtheorem{lemma}[theorem]{Lemma}

\theoremstyle{definition}
\newtheorem{definition}[theorem]{Definition}
\newtheorem{example}[theorem]{Example}
\newtheorem{remark}[theorem]{Remark}

\numberwithin{equation}{section}
\begin{document}
\title{Admissible perturbations of a multivalued Picard operator: \'Ciri\'c contraction condition; fixed point and stability results}


\author[]{\fnm{Cristina} \sur{Gheorghe}\email{cristina.gheorghe@ubbcluj.ro}

\affil[]{ \orgname{Faculty of Mathematics and Computer Science, Babeș-Bolyai
University}, \orgaddress{ \city{Cluj-Napoca}, \postcode{400084}, \country{Romania}}}

\pacs[2010 Mathematics Subject Classification]{47H10, 54H25.}}
\keywords{strict fixed point, multi-valued operator, admissible perturbation, convexity operator, data dependence, Ulam-Hyers stability, well-posedness, Ostrowski property, complete metric space, Picard operators, contraction-type conditions.}

\abstract{
This paper studies strict fixed point and stability results for multivalued operators which does not satisfy a \'Ciri\'c type contraction condition, but their admissible perturbation does. We focus on the conditions imposed on the admissible perturbation $T_G$ of a Picard operator $T:X\rightarrow P(X)$ such that the strict fixed point and stability results still hold for T. The results obtained are reformulated in terms of admissible perturbations in the sense of Takahashi and illustrated with some examples.
}

\maketitle


\section{Introduction}
\quad The concept of admissible perturbation of an operator was introduced by I.A. Rus in 2012, in \cite{Rus2012}, where iterative algorithms in terms of admissible perturbations of a single-valued operator were formulated and some data dependence and stability results were given in the context of these algorithms. Recently, in \cite{Rus},  some data dependence and stability results related to retraction-displacement conditions and admissible perturbations of a single-valued weakly Picard mapping were given.\\ 
\indent Lately, in \cite{art1}, some strict fixed point results related to multi-valued Picard operators were given, satisfying contraction conditions of \'Ciri\' c and \'Ciri\'c-Reich-Rus type.  The present work extends, for the case of admissible perturbations, the results related to multivalued Picard operators for whom a contraction condition of \'Ciri\' c type holds. The aim of this paper is to analyze the conditions that should be imposed on the admissible perturbation, so that the following conclusions hold true: the set of strict fixed points coincides with the set of fixed points and both have a unique element,  the sequence of Picard iterations converges to that unique strict fixed point, and the hypotheses under which a retraction-displacement condition occurs for the multivalued operator $T:X\rightarrow P(X)$. Then, based on the contraction principle formulated in terms of admissible perturbations, the data dependence and stability results are studied. We give similar results for the particular case of admissible perturbations in the sense of Takahashi and illustrate the results with relevant examples. \\
\indent First, let us introduce the main notions. 

\begin{definition} (see, e.g., \cite{Rus2012})\label{d0}
    Let X be a non-empty set and let $T:X\rightarrow P(X)$ represent a multivalued operator. Assume $G:X\times X \rightarrow X$ is an operator with the following properties:
\begin{itemize}
    \item [i)] $G(x,x)=x, $ for all $x\in X$;
    \item [ii)] $x,y\in X$ and $G(x,y)=x$ imply $y=x$.
\end{itemize}
Then a multivalued operator $T_G :X\rightarrow P(X)$ is defined as
\begin{equation}\label{tg}
  T_G(x)=G(x,T(x)):=\{G(x,u):u\in T(x)\}.  
\end{equation}
\end{definition}
\begin{definition}(see \cite{Rus2012})
    Let X be an non-empty space. Providing that the operator $G:X\times X\rightarrow X$ satisfies conditions i) and ii) from Definition \ref{d0},  the operator $T_G$ given in \eqref{tg} is known as the admissible perturbation of T corresponding to G.
\end{definition}
\begin{example} (see, e.g.,\cite{Petrusel})
    Let $(V,+,\mathbb{R})$ represent a vector space and let $X\subset V$ be a convex set. Suppose $\lambda\in (0,1)$ is a fixed parameter, let $T:X\rightarrow P(X)$ be a multivalued operator, and assume $G:X\times X\rightarrow X$ is an operator defined by $$G(x,y):=\lambda x+(1-\lambda) y.$$
    Then $T_G$ is an admissible perturbation of T corresponding to $G$.
\end{example}

\begin{lemma}\label{lema1} (see \cite{Berinde2016})
    Let X be a non-empty set and $T:X\rightarrow P(X)$ a multivalued operator. For an admissible perturbation $T_G$ corresponding to an operator $G:X\times X\rightarrow X$, the following affirmations hold true:
    \begin{itemize}
        \item[i)] $Fix(T_G)=Fix(T)$;
        \item[ii)] $SFix(T_G)=SFix(T)$.
    \end{itemize}
\end{lemma}
Let $(X,d)$ represent a complete metric space, suppose $T:X\rightarrow X$ is a multivalued operator and let $G(\cdot, T(\cdot))$ be an admissible perturbation of T. In this paper, our aim is to study the sufficient conditions imposed on the admissible perturbation $T_G$ such that the following properties hold true:
\begin{itemize}
    \item (DDE) data dependence estimate for the strict fixed point of T;
    \item (UH) Ulam-Hyers stability for the equation $T(x)=\{x\}$;
    \item (WP) well-posedness of the strict fixed point problem for T;
    \item (OP) Ostrowski property of the multivalued operator T;
    \item (QC) quasi-contraction condition for the multivalued operator T.
\end{itemize}
For other results related to admissible perturbations of a multivalued operator, we refer to \cite{ TakahashiArticle,articleGraphContr}.
\section*{Terminology and notations}
 Let $(X,d)$ be a metric space. Throughout this paper, $P(X)$ represents the family of all non-empty subsets of X, $P_{cp}$ defines the family of all non-empty compact subsets of X, and $P_{cl}$ is the family of all non-empty closed subsets of X. Also, for a multivalued operator $T:X\rightarrow P(X)$, the following notations are used :
\begin{itemize}
    \item the fixed point set of T:
    \begin{gather*}
         Fix(T):=\{x\in X: x\in T(x)\}
    \end{gather*}
    \item the strict fixed point set of T:
    \begin{gather*}
        SFix(T):=\{x\in X: \{x\}= T(x)\}
    \end{gather*}
    \item the graphic of T:
    \begin{gather*}
        Graph(T):=\{(x,y)\in X\times X: y\in T(x)\}
    \end{gather*}
\end{itemize}

We also recall, in the context of a metric space,  the definitions of the following functionals:
\begin{itemize}
    \item the gap functional generated by $d$:
\begin{gather*}
    D:P(X)\times P(X)\rightarrow \mathbb{R}_{+
},~D(A,B):=\inf\{d(a,b): a\in A, b\in B\}
\end{gather*}

\item the excess functional of A over B induced by d:
\begin{gather*}
    e:P(X)\times P(X)\rightarrow \mathbb{R}_{+}
\cup\{+\infty\}, ~e(A,B):=\sup\{D(a,B): a\in A \};
\end{gather*}

\item the Pompeiu-Hausdorff functional generated by $d$:
\begin{gather*}
    H:P(X)\times P(X)\rightarrow \mathbb{R}_{+
}\cup\{+\infty\},~ H(A,B):=\max\{e(A,B), e(B,A)\}.
\end{gather*}

\end{itemize}

  \begin{definition} (see \cite{AGPetruselSomeVariants})
    The operator $T:X\rightarrow P(X)$ is said to be  a multivalued Picard operator if the following statements hold true:
    \begin{itemize}
        \item[i)] $SFix(T)=Fix(T)=\{x^*\}$;
        \item[ii)] $T^n(x)\xrightarrow{H}\{x^*\}$ as $n\rightarrow \infty$, for each $x\in X.$ 
    \end{itemize}
\end{definition}
\begin{definition} (see \cite{Berinde2016})
  Let $(X,d)$ be a metric space, and let $T:X\rightarrow P(X)$ be a multivalued Picard operator. Then, the strong retraction-displacement condition is satisfied if there exists an increasing mapping $\Psi:\mathbb{R_+}\rightarrow \mathbb{R_+}$, continuous at 0 and $\Psi(0)=0$, such that the following condition holds true:
        \begin{equation*}\label{condRetrDispl}
            d(x,x^*)\leq \Psi(D(x,T(x)),\text{ for all } x\in X.
        \end{equation*}
    \end{definition} 
   We refer to \cite{Petrusel,TakahashiArticle,Carte55Rus, AGPetruselSomeVariants} for theorems related to multivalued Picard operators and \cite{Berinde2016, Rus} for results regarding the retraction-displacement condition.\\Let us recall some preliminary results related to multivalued operators that satisfy a given contraction-type condition. The interested reader can find them in detail in \cite{art1}.

\begin{theorem}\label{th55}(see \cite{art1})
    Let (X,d) be a complete metric space, and let $T:X\rightarrow P_{cl}(X)$ be a multivalued operator with $SFix(T)\neq \emptyset$. Assume there exist $\alpha, \beta, \gamma\geq 0$ with $\alpha+\beta+\gamma<1$ such that T satisfies the following condition:
    \begin{equation}\label{eqth55}
        H(T(x), T(y))\leq \alpha d(x,y)+\beta D(x, T(y))+\gamma D(y, T(x)), \text{ for all } x, y \in X.
    \end{equation}
   Then, the following relations hold:
    \begin{itemize}
        \item[i)] $Fix(T)=SFix(T)=\{x^*\}$;
        \item[ii)] For each $x\in X$, the sequence of sets $(T^n(x))_{n\in \mathbb{N}}$ converges to $\{x^*\}$ (with respect to the Pompeiu-Hausdorff generalized metric H);
        \item[iii)] For all $x\in X$ and for some $\xi \in (0,1)$, $$d(x,x^*)\leq \dfrac{1+\gamma}{(1-\alpha-\beta)\xi} D(x, T(x)).$$ 
    \end{itemize}
\end{theorem}
We mention that, in \cite{AGPetruselSomeVariants}, a similar strict fixed point principle was given for $\alpha-$ contractions. In the following, we give an example of an operator that is not a $\alpha-$ contraction but satisfies Theorem \ref{th55}.
\begin{example}
    Let $(X,d)$ be a complete metric space with $X=\Big[-\dfrac{8}{9}, \dfrac{8}{9}\Big]$ and let $T:X\rightarrow P_{cl}(X)$ be a multivalued operator, defined as $T(x)=[-x^2, x^2].$
    One easily sees that T is a \'Ciri\'c type contraction.\\ Additionally,  the following statements hold true:
    \begin{itemize}
        \item[i)] $Fix(T)=SFix(T)=\{0\}$;
        \item[ii)] $T^n(x)=[-x^{2n}, x^{2n}]$ and converges to $\{0\}$ with respect to Pompeiu-Hausdorff metric H, as $n\rightarrow \infty$; 
        \item[iii)] For some $\xi \in (0,1)$, one can get the retraction-displacement condition:
    \begin{align*}
        \begin{cases}
            |x|\leq\dfrac{1+\gamma}{(1-\alpha-\beta)\xi}|x-x^2|, &\text{for x} \in \Big[- \dfrac{8}{9},0\Big);\\
            ~\\
            |x|\leq\dfrac{1+\gamma}{(1-\alpha-\beta)\xi}|x+x^2|, &\text{for x} \in \Big[0, \dfrac{8}{9}\Big].
        \end{cases} 
    \end{align*} 
    \end{itemize} 
\end{example} 
As a consequence of Theorem \ref{th55}, we have the following result:
\begin{corollary}\label{lemmaProof}
    Let (X,d) be a complete metric space, and let $T:X\rightarrow P_{cl}(X)$ be a multivalued operator with $SFix(T)\neq \emptyset$. Suppose that there exist $\alpha, \beta, \gamma \geq 0 $ with $\alpha+\beta+\gamma<1$, satisfying the contraction-type condition (\ref{eqth55}). Then the following statements hold:
    \begin{itemize}
      \item [a)] $SFix(T)=Fix(T)=\{x^*\};$
        \item [b)] $H(T(x), \{x^*\})\leq \dfrac{\alpha+\beta}{1-\gamma} d(x,x^*),$ for all $x\in X.$
    \end{itemize}
\end{corollary}

\section{Admissible perturbation of a multivalued \'Ciri\'c type operator}
We aim to obtain strict fixed point and stability results for admissible perturbations presented in Definition \ref{d0}, using Theorem \ref{th55} and Lemma \ref{lema1}.
\begin{definition}
    Let $T:X\rightarrow P(X)$ a multivalued operator. Let $G:X\times X\rightarrow X$ and $T_G:X\rightarrow P(X)$ be an admissible perturbation of T corresponding to G. Then, for all $x\in X$,
    \begin{align*}
         T_{G}^{n}(x)=T_{G}(T_{G}^{n-1}(x))&=\bigcup_{y\in T_G^{n-1}(x)} T_G(y)\\&=\bigcup_{y\in T_G^{n-1}(x)}G(y, T(y)).
    \end{align*}
   
\end{definition}

\begin{lemma}\label{l3}
    Let $(X,d)$ represent a complete metric space, and let $G:X\times X\rightarrow X$ be a mapping such that the following conditions hold true:
    \begin{itemize}
    \item [i)] $G(x,x)=x, $ for all $x\in X$;
    \item [ii)] $x,y\in X$ and $G(x,y)=x$ imply $y=x$.
\end{itemize}
    Assume that $T:X\rightarrow P_{cl}(X)$  is a multivalued operator with $SFix(T)\neq \emptyset$ and consider $T_G$ is an admissible perturbation of T with respect to $T_G$. Provide there exists $l\in (0,1)$ such that 
    \begin{equation}\label{eql3}
    H(T(x), \{x^*\})\leq l H(T_G(x), x^*),\text{ for all } x\in X.
    \end{equation}
    Then, 
    \[H(T(Y), \{x^*\})\leq l H(T_G(Y), \{x^*\}), \text{ for all } Y\in P_{cl}(X).\]
    \begin{proof}
         Suppose $Y\in P_{cl}(X)$ and arbitrarily choose an element $v\in T(Y)$. Then there exists $y\in Y$ such that $v\in T(y).$ Hence, using \eqref{eql3}, we get
        \[d(v, x^*)\leq H(T(y), x^*)\leq l H(T_G(y), \{x^*\})\leq l H(T_G(Y), \{x^*\}).\]
        By taking the supremum of v over the set $T(Y)$, we obtain \[H(T(Y), \{x^*\})\leq l H(T_G(Y), \{x^*\}), \text{ for all } Y\in P_{cl}(X).\] \end{proof}
\end{lemma}
\begin{lemma}\label{l4}
    Let $(X,d)$ be a complete metric space and $G:X\times X\rightarrow X$ be a mapping satisfying the conditions:
    \begin{itemize}
    \item [i)] $G(x,x)=x, $ for all $x\in X$;
    \item [ii)] $x,y\in X$ and $G(x,y)=x$ imply $y=x$.
\end{itemize} 
Let $T:X\rightarrow P_{cl}(X)$ be a multivalued operator with $SFix(T)\neq \emptyset$ and let $T_G$ be an admissible perturbation of T with respect to G. Suppose that there exists $k\in (0,1)$ such that for all $x\in X$:
    \begin{equation}\label{eql4}
        H(T_G(x),\{x^*\})\leq k d(x,x^*).
    \end{equation}
    Then 
    \[H(T_G(Y), \{x^*\})\leq k H(Y, \{x^*\}), \text{ for all } Y\in P_{cl}(X).\]
    \begin{proof}
        Assume $Y\in P_{cl}(X)$ and let a $v\in T_G(Y)$ be an arbitrary element. Then, there exists $y\in Y$ such that $v\in T_G(y).$\\
       Moreover, for some $k\in (0,1)$, using \eqref{eql4}, we get
        \[d(v, x^*)\leq H(T_G(y), \{x^*\})\leq kd(y, \{x^*\})\leq H(Y, \{x^*\}). \]
        Taking the supremum of v over the set $T(Y)$, we obtain
        \[H(T_G(Y), \{x^*\})\leq k H(Y, \{x^*\}),\text{ for all } Y\in P(X).\]
    \end{proof}     
\end{lemma}
\begin{theorem}\label{th55tg}
    Let (X,d) be a complete metric space and $G:X\times X\rightarrow X$ be a mapping satisfying the conditions:
    \begin{itemize}
    \item [a)] $G(x,x)=x, $ for all $x\in X$;
    \item [b)] $x,y\in X$ and $G(x,y)=x$ imply $y=x$.
\end{itemize} Let $T:X\rightarrow P_{cl}(X)$ be a multivalued operator with $SFix(T)\neq\emptyset$ and $T_G$ be an admissible perturbation of T with respect to G. Suppose that there exist $\alpha, \beta,\gamma \geq 0 $ such that $\alpha+\beta+\gamma<1$ and
    \begin{equation}\label{eq2}
        H(T_G(x),T_G(y))\leq \alpha d(x,y)+\beta D(x, T_G(y))+\gamma D(y,T_G(x)), \text{ for all } x,y\in X.
    \end{equation} 
   Then the following statements hold true:
   \begin{itemize}
       \item[i)] $Fix(T)=SFix(T)=\{x^*\}$;\\
       \item[ii)] Additionally, assume there exists $l\in (0,1)$ such that 
       \begin{equation}\label{cond1ii)}
           H(T(x),\{x^*\})\leq l H(T_G(x), \{x^*\}), \text{ for all } x\in X.
       \end{equation}
        Then the sequence of successive approximations $(T^n(x))_{n\in \mathbb{N}}$ converges to $\{x^*\};$\\
       \item[iii)] Suppose there exists $L>0 $ such that the following inequality holds true for all $x\in X$:
       \begin{equation}\label{cond1iii)}
           D(x, T_G(x))\leq LD(x, T(x)).
       \end{equation}
       Then one has the retraction-displacement condition:
       \[d(x, x^*)\leq \frac{(1+\gamma)L}{(1-\alpha-\beta)\xi} D(x, T(x)), \text{ for all } x\in X \text{ and for some } \xi\in (0,1).\]
       
   \end{itemize}
   \begin{proof}
       From Lemma \ref{lema1} and Theorem \ref{th55}, affirmation i) follows immediately.\\
        ii) Let $x\in X$ be an arbitrary element.\\
    From Corollary \ref{lemmaProof} (applied for the multivalued $T_G$ having the unique strict fixed point $x^*$), we obtain that:
\begin{equation}\label{eq4}
    H(T_G(x),\{x^*\})\leq k d(x, x^*),
\end{equation}
where $k:=\dfrac{\alpha+\beta}{1-\gamma}\in (0,1).$\\
Let us consider the sequence $(T^n(x))_{n\in \mathbb{N}}$ induced by the multivalued operator T. We aim to prove by mathematical induction the proposition:
\[P(n):H(T^n(x), \{x^*\})\leq (lk)^n d(x, x^*), \text{ for all } n\in \mathbb{N}, n\geq 1.\]
From the hypothesis of ii) and from (\ref{eq4}), P(1) holds true. We suppose P(n) holds true and prove that $P(n)$ implies $P(n+1).$ Using Lemma \ref{l3} ii), Lemma \ref{l4} and condition \eqref{eq4}, we obtain:
\begin{align*}
    H(T^{n+1}(x),\{x^*\})&=H(T(T^n(x)), \{x^*\})\\&
  \leq l H(T_G(T^n(x)),\{x^*\})\\&
    \leq l\cdot k H(T^n(x), \{x^*\})\\&
    \overset{P(n)}\leq(lk)^{n+1} d(x,x^*).
    \end{align*}
Therefore, for all $x\in X$ and $n\in \mathbb{N}, n\geq 1,$
\[H(T^n(x),\{x^*\})\leq (lk)^n d(x, x^*),\]
which converges to 0 as $n\rightarrow \infty.$\\
iii) From Theorem \ref{th55}, for all $x\in X$, and some $\xi\in (0,1)$, we have
\begin{align*}
    d(x, x^*) &\leq \frac{1+\gamma}{(1-\alpha-\beta)\xi} D(x, T_G(x))\\&
    \leq \frac{1+\gamma}{(1-\alpha-\beta)\xi} LD(x, T(x))\\&
    = \frac{(1+\gamma)L}{(1-\alpha-\beta)\xi}D(x,T(x)).
\end{align*}
\end{proof}
\end{theorem}
The subsequent example is relevant for the circumstances specified in Theorem \ref{th55tg}.
\begin{example}
    Let $(X,d)$ be a complete metric space, with $X=\Bigg[\dfrac{1}{4}, 4\Bigg] $ and let $T:X\rightarrow P_{cl}(X)$ be a multivalued operator, such that 
  \begin{align*}
   T(x)=
      \begin{cases}
       \Bigg[1, \dfrac{1}{\sqrt{x}}\Bigg], &\text{ for } x\in \Bigg[\dfrac{1}{4},1\Bigg),\\
       ~\\
     [1,\sqrt{x}], &\text{ for } x\in [1,4].
       \end{cases}       
  \end{align*}
   One can easily observe that, by choosing $x=\dfrac{1}{4}$ and $y=1$ in Theorem \ref{th55}, the \'Ciri\'c contraction condition is not satisfied.\\ Let $T_G:X\rightarrow P_{cl}(X), T_G(x)=\dfrac{3}{4} x+\dfrac{1}{4} T(x)$ be an admissible perturbation of T.
   \begin{align*}
   T_G(x)=
       \begin{cases}
         \Bigg[\dfrac{3}{4}x+\dfrac{1}{4}, ~~~~~~\dfrac{3}{4}x+\dfrac{1}{4}\cdot\dfrac{1}{\sqrt{x}}\Bigg],  &\text{ for } x\in \Bigg[\dfrac{1}{4},1\Bigg),\\
       ~\\
       \Bigg[\dfrac{3}{4}x+\dfrac{1}{4},~~~~~~ \dfrac{3}{4}x+\dfrac{1}{4}\sqrt{x}\Bigg], &\text{ for } x\in [1,4].
       \end{cases}
   \end{align*}
For all $x, y\in X,$ there exist $\alpha, \beta, \gamma >0$  with $\alpha+\beta+\gamma<1$ such that $T_G$ satisfies the \'Ciri\'c contraction condition. Therefore, the upcoming affirmations hold true: 
\begin{itemize}
    \item[i)] $ SFix(T)=Fix(T)=\{1\}$;
    \item[ii)] There exists $l\in\Bigg(
    \dfrac{2}{3},1\Bigg)$ such that 
    \[H(T(x),\{1\})\leq l H(T_G(x), \{1\}), \text{ for all } x\in X,\]
    and the sequence of successive approximations $(T^n(x))_{n\in\mathbb{N}}$ converges to $\{1\};$
    \item[iii)] There exists $L\geq \dfrac{1}{4}$ such that, for all $x\in X,$
    \[D(x, T_G(x))\leq L D(x,T(x)),\]
    and one can obtain the retraction displacement condition for T:
       \[|x-1|\leq \dfrac{(1+\gamma)L}{(1-\alpha-\beta)\xi}|x-1|, \text{ for all } x\in X \text{ which is trivial for some } \xi \in (0,1). \]
\end{itemize}
\end{example}
In the following, we consider $T_G$ satisfying \'Ciri\'c contraction condition (\ref{eq2}) and analyze whether the fixed point and stability results suggested in Introduction hold. 
\subsection{Data dependence of the strict fixed point}
\begin{definition}\label{def DDE} (see \cite{article})
    Let $T:X\rightarrow P(X)$ a multivalued operator and let $F:X\rightarrow P(X)$ be another multivalued operator satisfying the following conditions:
    \begin{itemize}
        \item[i)]
            $SFix(F)\neq \emptyset$;\label{cond1}
        \item[ii)] there exists $\eta>0$ such that $H(T(x), F(x))\leq \eta,$ for all $x\in X.$
    \end{itemize}
    The strict fixed point problem $T(x)=\{x\}$ has the data dependence property if for each multivalued operator F satisfying conditions i) and ii) and for each strict fixed point $x^*$ of T, there exists a strict fixed point $u^*\in SFix(F)$ with the following property: \[\text{ there exists some } c> 0 \text{  such that } d(x^*, u^*)\leq c\eta.\]
\end{definition}
\begin{theorem}
    Let $(X,d)$ be a complete metric space and let $G:X \times X \rightarrow X$ be a mapping satisfying the conditions :
    \begin{itemize}
    \item [a)] $G(x,x)=x, $ for all $x\in X$;
    \item [b)] $x,y\in X$ and $G(x,y)=x$ imply $y=x$.
\end{itemize}
Let $T:X \rightarrow P_{cl}(X)$ be a multivalued operator with $SFix(T) \neq \emptyset $ and let $T_G$ be an admissible perturbation of T with respect to G. Supposing that there exist $\alpha, \beta, \gamma \geq 0$ such that $\alpha+\beta+\gamma <1 $ and
     \begin{align*}
        H(T_G(x),T_G(y))&\leq \alpha d(x,y)+\beta D(x, T_G(y))+\gamma D(y,T_G(x)), \text{ for all } x,y\in X,
    \end{align*} 
    and assuming that
    \[D(x,T_G(x))\leq L D(x, T(x)), \text{ for some } L>0,\]
   one has that the strict fixed point problem $T(x)=\{x\}$ has the data dependence property. 
    \begin{proof}
       From Theorem \ref{th55tg}, it follows directly that $SFix(T)=\{x^*\}$ and 
        \[d(x, x^*)\leq \widetilde{K} D(x,T(x)),\]
        where $\widetilde{K}:=\dfrac{1+\gamma}{(1-\alpha-\beta)\nu},$ $\widetilde{K}>0$ and $\nu\in (0,1).$\\
        Let $F:X\rightarrow P(X)$ be a multivalued operator satisfying i) and ii) from Definition \ref{def DDE}.\\ Suppose $u^*\in SFix(F)$. From relation ii) of Definition \ref{def DDE},  we infer 
        \begin{align*}
            d(u^*, x^*)&\leq \widetilde{K} D(u^*, T(u^*))\\&
            =\widetilde{K} H(F(u^*), T(u^*)){\leq} \widetilde{K} \eta.
        \end{align*}
    \end{proof}
\end{theorem}
In the following, we introduce a more general result for data dependence estimate, without using Theorem 9. A similar result was obtained in \cite{Rus} for the univalued case. We shall introduce the background. \\
\indent Let $(X, d)$ a metric space, $T, F:X\rightarrow P(X)$ be multivalued operators with $SFix(T)=\{x^*\} $ and $SFix(F)\neq \emptyset$. We aim to study which conditions should be imposed on the admissible perturbation of T, ensuring there exists an increasing mapping $\Psi:\mathbb{R}_{+}\rightarrow\mathbb{R}_{+}$, continuous at zero with $\theta(0)=0$, such that
\[d(y^*,x^*)\leq \Psi(\eta), \text{ for all } y^*\in SFix(F).\]
\begin{theorem}\label{ThCondPicard}
     Let (X,d) be a complete metric space,  $T:X\rightarrow P(X)$ a multivalued operator with $SFix(T)=\{x^*\}$ and $T_G$ an admissible perturbation.\\ 
    Let $F:X\rightarrow P(X)$ be a multivalued operator such that $SFix(F)\neq \emptyset$ and
\[H(T(x),F(x))\leq \eta, \text{ for all 
} x\in X, \text{ for some } \eta\in \mathbb{R}_{+}^{*}.\]
We suppose that:
    \begin{itemize}
        \item[1)] $T_G$ is a $\Psi-MP$ operator with $SFix(T_G)=\{x^*\}$;
        \item[2)] $D(x, T_G(x))\leq c D(x, T(x)), \text{ for all } x\in X \text{ with some } c\in \mathbb{R}_{+}^{*};$
        \item[3)] there exists $\eta > 0$ such that $H(T(x),F(x))\leq \eta,\text{ for all } x\in X$.
        \end{itemize}
Then the following statements are true:
\begin{itemize}
    \item[i)] $d(x,x^*)\leq \Psi(c D(x,T(x))
    ), \text{ for all } x\in X; $
    \item[ii)] $d(x^*,y^*)\leq \Psi (c\eta) \text{ for all } y^*\in SFix(F).$
\end{itemize}
    \begin{proof}
        i) From Lemma \ref{lema1}, we deduce
        \[SFix(T_G)=SFix(T)=\{x^*\}.\]
        Since $T_G$ is a $\Psi-$ MP operator, we infer
        \[d(x,x^*)\leq \Psi (D(x,T_G(x))), \text{ for all } x\in X.\]
        Relation 2) now implies  
        \[d(x,x^*)\leq \Psi(cD(x, T(x))), \text{ for all } x\in X.\]
        ii) Let $x=y^*\in SFix(F)$. From i), we get
        \begin{align*}
             d(y^*,x^*)&\leq \Psi (cD(y^*,T(y^*)))=\Psi(cH(F(y^*), T(y^*)))\\&\leq \Psi (c\eta), 
        \end{align*}
         for all $x\in X$  and for some $\eta\in \mathbb{R^*_{+}}.$
    \end{proof}
\end{theorem}

\subsection{Ulam-Hyers stability property of the strict fixed point problem}
\begin{definition}\label{def UH} (see \cite{Carte55Rus})
    The strict fixed point problem $T(x)=\{x\}$ is said to be Ulam-Hyers stable if there exists $c\in \mathbb{R^*}$ such that for any solution $y^*$ of the given equation
    $$D(y, T(y))\leq \varepsilon,$$
    there exists a strict fixed point of T, namely $x^*\in SFix(T)$, satisfying the inequality
    $$d(y^*, x^*)\leq \varepsilon.$$
\end{definition}
\begin{theorem}
    Let $(X,d)$ be a metric space and let $T:X\rightarrow P(X)$ be a multivalued operator with $SFix(T)\neq \emptyset$. Suppose that $G:X\times X\rightarrow  X$ is a mapping satisfying the conditions:
    \begin{itemize}
    \item [i)] $G(x,x)=x, $ for all $x\in X$;
    \item [ii)] $x,y\in X$ and $G(x,y)=x$ imply $y=x$,
\end{itemize}
Consider $T_G$ as the admissible perturbation of T with respect to G. 
    We denote, for some $\varepsilon>0$, by $y^*$ the solution of the inequation
    \[D(y, T(y))\leq \varepsilon.\]
    Assume $T_G$ satisfies the hypotheses of Theorem \ref{th55tg} and  the following additional condition holds:
    \[\text{iii) there exists } L>0 \text{ such that } D(x, T_G(x))\leq L D(x, T(x)), \text{ for every } L>0. \]
    Then the strict fixed point problem  $T(x)=\{x\}$ is Ulam-Hyers stable.
    \begin{proof}
        Let $y^*$ be a solution of the inequation $D(y, T_G(y))\leq \varepsilon.$ From Theorem \ref{th55}, we have that $SFix(T_G)=\{x^*\}.$ Then,
        \begin{align*}
            d(y^*, x^*)&=D(y^*, T_G(x^*))
            \\&\leq D(y^*, T_G(y^*))+H(T_G(y^*), T_G(x^*))\\&
            \leq LD(y^*, T(y^*))+\alpha d(x^*, y^*)+\beta D(x^*, T_G(y^*))+\gamma D(y^*, T_G(x^*))\\&
            \leq L\varepsilon +(\alpha+\beta+\gamma)d(x^*, y^*)+\beta LD(y^*, T(y^*))\\&
            \leq L\varepsilon (1+\beta)+(\alpha+\beta+\gamma) d(x^*, y^*).
        \end{align*}
        Therefore, 
\begin{equation*}
    d(x^*, y^*)\leq c\varepsilon,
\end{equation*}
         where 
         $$c:=\dfrac{L(1+\beta)}{1-\alpha-\beta-\gamma}\in \mathbb{R^*_+}.$$
     \end{proof}
\end{theorem} 
\subsection{Well-posedness of the strict fixed point problem}
\begin{definition}\label{def WP} (see, e.g., \cite{article}, \cite{artprec})
    The strict fixed point problem $T(x)=\{x\}$ is well-posed in the sense of Reich and Zaslavski (see \cite{Reich2005ANO}) if T has a unique strict fixed point $x^*$ and the following implication holds true: \[(u_n)_{n\in\mathbb{N}}\subset X \text{ and } D(u_n, T(u_n))\rightarrow 0\text{ implies } u_n\rightarrow x^* \text{ as } n\rightarrow \infty.\]
\end{definition}
\begin{theorem}
    Let (X,d) be a metric space and let $T:X\rightarrow P_{cl}(X)$ be a multivalued operator with $SFix(T)\neq \emptyset$. Suppose that $G:X\times X\rightarrow  X$ is a mapping that satisfies the following conditions:
    \begin{itemize}
    \item [a)] $G(x,x)=x, $ for all $x\in X$;
    \item [b)] $x,y\in X$ and $G(x,y)=x$ imply $y=x$.
\end{itemize}
Consider $T_G$ as the admissible perturbation of T with respect to G.
    Assume $T_G$ satisfies the hypotheses of Theorem \ref{th55tg}, including the additional conditions ii) and iii), namely:
     \begin{itemize}
      
       \item[ii)] there exists $l\in (0,1)$ such that 
       \begin{equation}\label{cond2ii)}
           H(T(x),\{x^*\})\leq l H(T_G(x), \{x^*\}), \text{ for all } x\in X;
       \end{equation}

       \item[iii)] there exists $L>0 $ such that, for all $x\in X$, the following inequality holds:
       \begin{equation}\label{cond2iii)}
           D(x, T_G(x))\leq LD(x, T(x)).
       \end{equation}
       
   \end{itemize}
    Then the strict fixed point problem $T(x)=\{x\}$ is well-posed.
    \begin{proof}
      From Theorem \ref{th55tg}, we have $SFix(T)=\{x^*\}$. Moreover, from condition iii), for all $x\in X,$ there exists $\eta\in (0,1)$ satisfying
      \[d(u_n, x^*)\leq \dfrac{1+\gamma}{(1-\alpha-\beta)\eta}D(u_n, T(u_n))\rightarrow 0 \text{ as } n\rightarrow \infty.\]
      Alternatively, we infer
      \begin{align*}
           d(u_n, x^*)&\leq D(u_n, T_G(u_n))+H(T_G(u_n),T_G(x^*))\\&
           \leq D(u_n, T_G(u_n))+\alpha d(u_n, x^*)+\beta [d(x^*, u_n)+D(u_n, T_G(u_n))]+\gamma D(u_n, x^*).
      \end{align*}
       Consequently, 
       \begin{align*}
           d(u_n, x^*)&\leq \dfrac{1+\beta}{1-\alpha-\beta-\gamma}D(u_n, T_G(u_n))\\&\leq L\dfrac{1+\beta}{1-\alpha-\beta-\gamma} D(u_n, T(u_n))\rightarrow 0\text{ as } n\rightarrow \infty.
       \end{align*}
    \end{proof}
\end{theorem}
\subsection{Ostrowski stability property of strict fixed point problem}
\begin{definition}\label{def OP} (see \cite{art9},\cite{Rus})
    The fixed point problem $T(x)=\{x\}$ is said to have the Ostrowski stability property if T has a unique strict fixed point $x^*$ and the following implication holds:\[(v_n)_{n\in\mathbb{N}}\subset X, D(v_{n+1}, T(v_n))\rightarrow 0\Rightarrow v_n\rightarrow x^* \text{ as } n\rightarrow \infty.\]
\end{definition}
\begin{lemma}\label{CauchyToeplitzlema}
    (Cauchy-Toeplitz lemma; see, e.g., \cite{CauchyToeplitz}) We consider $(a_n)_{n\in\mathbb{N}}$ a sequence of positive real numbers such that the series $\displaystyle\sum_{n=0}^\infty$ is convergent. Let $(b_n)_{n\in \mathbb{N}}$ be another sequence of nonnegative numbers convergent to zero. Then, we have:\[\lim_{n\rightarrow \infty} \Bigg(\sum_{k=0}^n a_{n-k} b_k\Bigg)=0.\]
\end{lemma}
\begin{theorem}
    Let (X, d) represent a complete metric space and let $T:X\rightarrow P(X)$ be a multivalued operator with $SFix(T)\neq \emptyset$. Assume that $G:X \times X \rightarrow  X$ is a mapping satisfying the conditions:
    \begin{itemize}
    \item [a)] $G(x,x)=x, $ for all $x\in X$;
    \item [b)] $x,y\in X$ and $G(x,y)=x$ imply $y=x$.
\end{itemize}
Suppose $T_G$ is an admissible perturbation of T with respect to G. We assume true the hypotheses of Theorem \ref{th55tg} related to $T_G$, including the supplementary conditions ii) and iii), namely:
\begin{itemize}
 \item[ii)] there exists $l\in (0,1)$ such that 
       \begin{equation}\label{cond3ii)}
           H(T(x),\{x^*\})\leq l H(T_G(x), \{x^*\}), \text{ for all } x\in X;
       \end{equation}

       \item[iii)] there exists $L>0 $ such that, for all $x\in X$, the following inequality holds:
       \begin{equation}\label{cond3iii)}
           D(x, T_G(x))\leq LD(x, T(x)).
       \end{equation}
   \end{itemize}
    Then the strict fixed point problem $T(x)=\{x\}$ is Ostrowski stable. 
    \begin{proof}
        We deduce, from Theorem \ref{th55} applied to $T_G$, that $SFix(T_G)=\{x^*\}.$ Successively, we get
       \begin{align*}
           d(v_{n+1}, x^*)&\leq D(v_{n+1}, T_G(v_n))+H(T_G(v_n), T_G(x^*))\\&
           \leq D(v_{n+1},T_G(v_n))+\alpha d(v_n, x^*)+\beta D(v_n, T_G(x^*))+\gamma D(T_G(v_n), x^*))\\&
           \leq D(v_{n+1},T_G(v_n))+(\alpha+\beta)d(v_n, x^*)+\gamma [d(v_{n+1}, x^*)+D(v_{n+1}, T_G(v_n))].
           \end{align*}
       With the notation $$k:=\dfrac{\alpha+\beta}{1-\alpha}, \text{ where } 0<k<1;$$ and using condition iii) from Theorem \ref{th55tg}, we obtain
       \begin{align*}
           d(v_{n+1},x^*)&\leq \frac{1+\gamma}{1-\gamma}D(v_{n+1}, T_G(v_n))+kd(v_n, x^*)\\&
           \leq \frac{1+\gamma}{1-\gamma} D(v_{n+1}, T_G(v_n))+k\Big[\frac{1+\gamma}{1-\gamma}D(v_n, T_G(v_{n-1}))+k d(v_{n-1}, x^*)\Big]\\&
          \leq  \cdots\\&
          \leq \frac{1+\gamma}{1-\gamma}\Big[D(v_{n+1}, T_G(v_n))+k D(v_n, T_G(v_{n-1}))+\cdots + k^n D(v_1, T_G(v_0))\Big]\\&
         \leq L \frac{1+\gamma}{1-\gamma}\Big[D(v_{n+1}, T(v_n))+\\&+k D(v_n, T(v_{n-1}))+\cdots +k^n D(v_1, T(v_0))\Big]. 
       \end{align*}
       From Cauchy-Toeplitz Lemma, having $\displaystyle \sum_{n=0}^\infty k^n$ convergent and $\displaystyle \lim_{n\rightarrow \infty} D(v_{n+1}, T(v_n))=0$, one deduces that $d(v_{n+1}, x^*)\rightarrow 0$ and therefore, the strict fixed point problem $T(x)=x$ is Ostrowski stable.
    \end{proof}
\end{theorem}
\subsection{Quasi-contraction condition for the multivalued operator T} 
\begin{definition}\label{def QC}
    Let $(X,d)$ be a complete metric space. We call quasi-contraction an operator $T:X\rightarrow P(X)$ with $SFix(T)\neq \emptyset$  and for whom there exists $l\in(0,1)$ satisfying
    \[H(T(x), x^*)\leq ld(x,x^*)\text{ for all } x^*\in SFix(T) \text{ and for all } x\in X.  \]
\end{definition}
\begin{theorem}\label{thstrongQC}
    Let $(X,d)$ represent a complete metric space, let $T:X\rightarrow P(X)$ be a multivalued operator with $SFix(T)\neq \emptyset$ and consider $G:X\times X\rightarrow  X$ be a mapping satisfying the conditions:
    \begin{itemize}
    \item [a)] $G(x,x)=x, $ for all $x\in X$;
    \item [b)] $x,y\in X$ and $G(x,y)=x$ imply $y=x$.
    \end{itemize}
    Suppose $T_G$ is an admissible perturbation of T with respect to G. Assume $T_G$ satisfies the hypotheses of Theorem \ref{th55tg} and the following supplementary condition holds:
       \begin{equation}\label{cond4ii)}
           \text{there exists $l\in (0,1)$  such that } H(T(x),\{x^*\})\leq l H(T_G(x), \{x^*\}), \text{ for all } x\in X.
       \end{equation}
  Then the multivalued operator T is a quasi-contraction.
    \begin{proof}
        From Theorem $\ref{th55}$, we deduce $SFix(T_G)=\{x^*\}$. One has:
        \begin{align*}
            D(T_G(x), x^*)\leq H(T_G(x), \{x^*\})& \leq \alpha d(x, x^*)+\beta D(x, T_G(x^*))+\gamma D(x^*, T_G(x))\\
            &\leq \alpha d(x, x^*)+\beta D(x, T_G(x^*))+\gamma H(x^*, T_G(x))
        \end{align*}
        Therefore, $$(1-\gamma)H(T_G(x), x^*)\leq (\alpha+\beta) d(x, x^*).$$ From  condition (\ref{cond4ii)}), there exists $l\in (0,1 )$ such that
        \begin{align*}
            H(T(x), \{x^*\})\leq l H(T_G(x), \{x^*\})\leq l \dfrac{\alpha+\beta}{1-\gamma} d(x, x^*). 
        \end{align*}
        Since  $l \dfrac{\alpha+\beta}{1-\gamma}\in (0,1)$, it follows that T is a quasi-contraction.
    \end{proof}
\end{theorem}
\begin{definition}\label{def weakQC} 
    Let $(X,d)$ be a complete metric space. We call weak quasi-contraction an operator $T:X\rightarrow P(X)$ such that $SFix(T)\neq \emptyset$  and for whom exists $l\in(0,1)$ satisfying
    \[D(T(x), x^*)\leq ld(x,x^*)\text{ for all } x^*\in SFix(T) \text{ and for all } x\in X.  \]
\end{definition}
\begin{remark}
 By substituting condition \eqref{cond4ii)} from Theorem \ref{thstrongQC} with
 \begin{equation}\label{condii)2}
       \text{there exists $l\in (0,1)$  such that } D(T(x),\{x^*\})\leq l D(T_G(x), \{x^*\}), \text{ for all } x\in X,
 \end{equation}
 we obtain that T is a weak quasi-contraction.
  \begin{proof}
           Similarly to the proof of Theorem \ref{thstrongQC}, we have $$(1-\gamma)D(T_G(x), x^*)\leq (\alpha+\beta) d(x, x^*).$$ From  condition (\ref{condii)2}), there exists $l\in (0,1 )$ such that
        \begin{align*}
            D(T(x), x^*)\leq l D(T_G(x), \{x^*\})\leq l \dfrac{\alpha+\beta}{1-\gamma} d(x, x^*). 
        \end{align*}
        Since  $l \dfrac{\alpha+\beta}{1-\gamma}\in (0,1)$ it follows that T is a weak quasi-contraction.
         \end{proof}
\end{remark}
The discussion can be broadened to include other contraction-type conditions.
\begin{remark}
        By replacing the contraction-type condition (\ref{eqth55}) with the \'Ciri\'c-Reich-Rus contraction condition:
        \begin{equation}\label{cond2}
            H(T(x), T(y))\leq \alpha d(x,y)+\beta D(x, T(x))+\gamma D(y, T(y)), \text{ for all } x, y\in X,     
        \end{equation}
         where $\alpha, \beta, \gamma >0$  and $\alpha+2\beta<1$, one can  obtain similar results for admissible perturbations of a multivalued operator. For some main results related to multivalued operators satisfying  the \'Ciri\'c-Reich-Rus contraction condition, we refer to \cite{art1}.
    \end{remark}
    \begin{remark}
        By joining together the contraction-type conditions (\ref{th55}) and (\ref{cond2}), one can obtain similar results for admissible perturbations of multivalued operators $T:X\rightarrow P(X)$, with $SFix(T)\neq \emptyset$, satisfying the contraction-type condition:
        \begin{equation}
            H(T(x, T(y))\leq \alpha d(x,y)+\beta[D(x, T(x))+D(y,T(y))]+\gamma [D(x, T(y))+D(y, T(x))],
             \end{equation}
            where $\alpha,\beta, \gamma>0$, with $\alpha+2\beta<1.$
       
    \end{remark}
     
    \begin{remark}
        Letting $\alpha\in (0,1)$ and $\beta=\gamma=0$ the discussion can be moved to admissible perturbations of an $\alpha$-contraction (see, e.g., Theorem 5.5 from \cite{AGPetruselSomeVariants}).
    \end{remark}
\section{Application: fixed point and stability results for an admissible perturbation in the sense of Takahashi}
\subsection{Main results for Takahashi admissible perturbations}
\begin{definition}
    Let $(X,d)$ be a metric space. Assume that the operator $W:X\times X\times [0,1]\rightarrow P_{cl}(X)$ satisfies the following property:
\[d(u,W(x,y,\lambda))\leq \lambda d(u,x)+(1-\lambda)d(u,y).\] 
Then the triplet $(X,d,W)$ is called a convex metric space in the sense of Takahashi.
\end{definition}
\begin{definition}\label{dW}
    Let X be a non-empty set and consider $T:X\rightarrow P_{cl}(X)$ a multivalued operator. Assume $W:X\times X\times [0,1]\rightarrow P_{cl}(X)$ such that $(X,d,W)$ is a convex metric space in the sense of Takahashi, having the following properties:
    \begin{itemize}
        \item[i)] $W(x,x,\lambda)=x,$ for all $x \in X$, $\lambda\in (0,1)$;
        \item[ii)] 
       
             $\lambda\in (0,1), x,y\in X \text{ and } W(x,y, \lambda)=x \text{ imply } y=x.$
    \end{itemize}
    The operator $T_W:X\rightarrow P_{cl}(X)$, defined as
    \[T_W(x)=\{W(x,y,\lambda):y\in T(x)\},\]
    is called the Takahashi admissible perturbation of T, corresponding to W.
\end{definition}
\begin{remark}
    We observe that $\displaystyle W(x, T(x), \lambda)=\bigcup_{y\in T(x)} W(x,y,\lambda)$, for all $x\in X$.
\end{remark}
\begin{theorem}\label{thTW}
    Let $(X,d)$ be a complete metric space and $W:X\times X\times [0,1]\rightarrow P_{cl}(X)$  a mapping satisfying:
    \begin{itemize}
    \item[i)] $W(x,x,\lambda)=x,$ for all $x \in X$, $\lambda\in (0,1)$;
        \item[ii)] 
             $\lambda\in (0,1), x,y\in X \text{ and } W(x,y, \lambda)=x \text{ imply } y=x.$ 
    \end{itemize}
    Let $T:X\rightarrow P_{cl}(X)$ be a multivalued operator with $SFix(T)\neq \emptyset$. Assume that $T_W$ is an admissible perturbation of T with respect to W.  Suppose there exist $\alpha, \beta, \gamma \geq 0$ such that $\alpha+\beta+\gamma<1$ and 
\begin{align*}
    H(W(x,T(x),\lambda), W(y, T(y), \lambda))&\leq \alpha d(x,y)+\beta D(x, W(y, T(y), \lambda))+\gamma D(y, W(x, T(x), \lambda)), \\&\text{ for all } x, y\in X.
\end{align*}
Then the following statements hold true:
\begin{itemize}
    \item[i)] $Fix(T)=SFix(T)=\{x^*\};$
    \item[ii)] Moreover, assume there exists $l\in (0,1)$ such that 
    \[H(T(x), \{x^*\})\leq l H(W(x, T(x), \lambda),\{x^*\}),\]
    for all $x\in X$ and $x^*\in SFix(T)$.\\
    Then the sequence of successive approximations $(T^n(x))$ converges to $x^*$;
    \item[iii)] Additionally, supposing there exists $L>0$ such that for all $x\in X,$
    \[D(x, W(x,T(x),\lambda))\leq L D(x, T(x)),\]
    then
    \[d(x, x^*)\leq \dfrac{(1+\gamma)L}{(1-\alpha-\gamma)\xi}D(x, T(x)),\]
    for all $x \in X$ and for some $\xi\in (0,1).$
\end{itemize}
    \begin{proof}
        The proof is trivial by applying Theorem \ref{th55tg} to $T_W$.
    \end{proof}

\end{theorem}
\begin{theorem}
    Let $(X,d)$ be a complete metric space and $W:X\times X \times [0,1] \rightarrow P_{cl}(X)$ be a mapping that meets the conditions: 
    \begin{itemize}
    \item[a)] $W(x,x,\lambda)=x,$ for all $x \in X$, $\lambda\in (0,1)$;
        \item[b)] 
             $\lambda\in (0,1), x,y\in X \text{ and } W(x,y, \lambda)=x \text{ imply } y=x.$ 
    \end{itemize}
    Let $T:X\rightarrow P_{cl}(X)$ be a multivalued operator with $SFix(T)\neq \emptyset$. Suppose $T_W$ is an admissible perturbation of T with respect to W. Provide there exists $\alpha, \beta, \gamma\geq 0$ such that $\alpha+\beta+\gamma<1$ satisfying 
    \begin{align*}
        H(W(x, T(x), \lambda)), W(y, T(y), \lambda))&\leq \alpha d(x,y)+\beta D(x, W(y, T(y), \lambda))+\gamma D(y, W(x, T(x), \lambda)) \\&\text{ for all } x, y\in X,
     \end{align*}
     and 
     \[D(x, W(x, T(x), \lambda))\leq L D(x, T(x)), \text{ for some } L>1.\]
     Then the strict fixed point problem $T(x)=\{x\}$ has the data dependence property. 
\end{theorem}
We can also formulate a more general result for the Takahashi perturbation:
\begin{theorem}
    Let $(X,d)$ be a complete metric space, $T,F:X\rightarrow P_{cl}(X)$ be multivalued operators such that $SFix (T)=\{x^*\}$ and $SFix(F)\neq \emptyset$ and 
    \[H(T(x), F(x))\leq \eta, \text{ for all } x\in X, \text{ for some } \eta \in \mathbb{R}_+^*.\]
    Assume $W:X\times X
    \times [0,1]\rightarrow X$ is a mapping satisfying the following conditions:
    \begin{itemize}
    \item[i)] $W(x,x,\lambda)=x,$ for all $x \in X$, $\lambda\in (0,1)$;
        \item[ii)] 
             $\lambda\in (0,1), x,y\in X \text{ and } W(x,y, \lambda)=x \text{ imply } y=x.$ 
    \end{itemize}
    Let $T_W$ representing an admissible perturbation of T in the sense of Takahashi.\\
    We suppose that:
    \begin{itemize}
        \item[1)] $T_W$ is a $\Psi-$ MP operator with $SFix(T_W)=\{x^*\};$
        \item[2)] $D(x, T_W(x))\leq c D(x, T(x)),$ for all $x\in X$ and for some $c\in \mathbb{R}_+^*$;
        \item[3)] there exists $\eta>0$ such that $H(T(x), F(x))\leq \eta$, for all $x\in X$.
    \end{itemize}
        Then, the following affirmations hold true:
        \begin{itemize}
            \item[i)] $d(x, x^*)\leq \Psi(cD(x,T(x)))$, for all $x\in X$;
            \item[ii)] $d(x^*, y^*)\leq \Psi(c \eta )$, for all $y^*\in SFix(F)$.
    \end{itemize}
\end{theorem}
\begin{theorem}
    Consider $(X,d)$ be a metric space, let $T:X\rightarrow P_{cl}(X)$ represent a multivalued operator with $SFix(T)\neq \emptyset$ and suppose $W:X\times X\times [0,1]\rightarrow P_{cl}(X)$ is a mapping satisfying the conditions: \begin{itemize}
    \item[a)] $W(x,x,\lambda)=x,$ for all $x \in X$, $\lambda\in (0,1)$;
        \item[b)] 
             $\lambda\in (0,1), x,y\in X \text{ and } W(x,y, \lambda)=x \text{ imply } y=x.$ 
    \end{itemize} 
    Let $T:X\rightarrow P_{cl}(X)$ be a multivalued operator with $SFix(T)\neq \emptyset$ and assume that $T_W$ is the admissible perturbation with respect to T. For some $\varepsilon>0$, we consider $y^*$ being the solution of the inequality
    \[D(y, T(y))\leq \varepsilon.\]
    Assume $T_W$ satisfies the hypotheses of Theorem \ref{thTW} and the following condition holds:
    \[D(x, T_W(x))\leq LD(x,T(x)), \text{ for some } L>1.\]
     Then the strict fixed point problem $T(x)=\{x\}$ is Ulam-Hyers stable.
\end{theorem}
\begin{theorem}
    Let $(X,d)$ be a metric space, consider $T:X\rightarrow P_{cl}(X)$ is a multivalued operator such that $SFix (T)\neq \emptyset.$ Suppose that $W:X\times X\times [0,1]\rightarrow P_{cl}(X)$ is a mapping satisfying the following conditions:
    \begin{itemize}
    \item[a)] $W(x,x,\lambda)=x,$ for all $x \in X$, $\lambda\in (0,1)$;
        \item[b)] 
             $\lambda\in (0,1), x,y\in X \text{ and } W(x,y, \lambda)=x \text{ imply } y=x;$ 
    \end{itemize}
    Let $T_W$ be an admissible perturbation of T with respect to W. Suppose $T_W$ satisfies the hypotheses of Theorem \ref{thTW}, including the additional conditions: 
    \begin{itemize}
    \item[ii)] there exists $l\in (0,1)$ such that 
    \[H(T(x), \{x^*\})\leq l H(W(x, T(x), \lambda),\{x^*\}),\]
    for all $x\in X$ and $x^*\in SFix(T)$;
    \item[iii)] there exists $L>0$ such that for all $x\in X,$
    \[D(x, W(x,T(x),\lambda))\leq L D(x, T(x)).\]   
\end{itemize}
    Then the fixed point problem $T(x)=\{x\}$ is well-posed. 
    \end{theorem}
    \begin{theorem}
        Let $(X,d)$ be a complete metric space and let $T:X\rightarrow P_{cl}(X) $ represent a multivalued operator with $SFix(T)\neq \emptyset$. Consider $W:X\times X\times [0,1]\rightarrow P_{cl}(X)$ is a mapping satisfying the conditions:
        \begin{itemize}
    \item[a)] $W(x,x,\lambda)=x,$ for all $x \in X$, $\lambda\in (0,1)$;
        \item[b)] 
             $\lambda\in (0,1), x,y\in X \text{ and } W(x,y, \lambda)=x \text{ imply } y=x;$ 
    \end{itemize}
        Let $T_W$ be an admissible perturbation of T, in the sense of Takahashi,  with respect to W. Assume the hypotheses of Theorem \ref{thTW} (applied to $T_W$) hold true, including the supplementary conditions from ii) and iii), namely: 
 \begin{itemize}
    \item[ii)] there exists $l\in (0,1)$ such that 
    \[H(T(x), \{x^*\})\leq l H(W(x, T(x), \lambda),\{x^*\}),\]
    for all $x\in X$ and $x^*\in SFix(T)$;
    \item[iii)] there exists $L>0$ such that for all $x\in X,$
    \[D(x, W(x,T(x),\lambda))\leq L D(x, T(x)).\]
    \end{itemize}
Then  the strict fixed point problem $T(x)=\{x\}$ has the Ostrowski stability property. 
    \end{theorem}
    \begin{theorem}
        Let $(X,d)$ represent a complete metric space. Assume $T:X\rightarrow P_{cl}(X)$ is a multivalued operator with $SFix(T)\neq \emptyset$ and let $W:X\times X\times [0,1]\rightarrow P_{cl}(X)$ be a mapping for which the following conditions hold true:
         \begin{itemize}
    \item[a)] $W(x,x,\lambda)=x,$ for all $x \in X$, $\lambda\in (0,1)$;
        \item[b)] 
             $\lambda\in (0,1), x,y\in X \text{ and } W(x,y, \lambda)=x \text{ imply } y=x.$ 
    \end{itemize}
        Consider $T_W$ is an admissible perturbation of T, in the sense of Takahashi, with respect to W. \\Provide the hypotheses of Theorem \ref{thTW} (applied to $T_W$) hold true, including the supplementary condition:
        \begin{itemize}
            \item[iii)] there exists $L>0$ such that for all $x\in X,$
    \[D(x, W(x,T(x),\lambda))\leq L D(x, T(x)).\]
        \end{itemize}
    Then the multivalued operator T is a quasi-contraction. 
    \end{theorem}
    \subsection{Example: the convexity operator in the sense of Takahashi}
    In the following, we will illustrate a concrete example of the presented theory. Let us consider the convexity operator defined by Takahashi in \cite{TakahashiW},  $$W:X\times X\times [0,1] \rightarrow P_{cl}(X),$$ with $$W(x,T(x), \lambda)=\lambda x+(1-\lambda) T(x).$$
    For results concerning convex structures, we refer to \cite{ArtGudder} and {TakahashiArticle}.
     \begin{remark}
         We note that $T_W(x) =\lambda x+(1-\lambda)T(x).$
         \begin{proof}
For any $x\in X$, we have
         \begin{align*}
             T_W(x) =W(x, T(x), \lambda)&=\displaystyle \bigcup_{y\in T(x)} W(x,y,\lambda)=\bigcup_{y\in T(x)} (\lambda x+(1-\lambda)y)\\&=\lambda x+(1-\lambda)\bigcup_{y\in T(x)} y\\&=\lambda x+(1-\lambda)T(x).
         \end{align*}
         \end{proof}
     \end{remark}
     \begin{example}
        Let $(X,d)$ be a complete metric space, where $X=\Bigg[\dfrac{1}{4}, 4\Bigg]$, and let $T:X\rightarrow P_{cl}(X)$ be a multivalued operator, defined as
        \begin{align*}
   T(x)=
      \begin{cases}
       \Bigg[1, \dfrac{1}{\sqrt{x}}\Bigg], &\text{ for } x\in \Bigg[\dfrac{1}{4},1\Bigg),\\
       ~\\
     [1, \sqrt{x}], &\text{ for } x\in [1,4].
       \end{cases}       
  \end{align*}
  Let $T_W:X\times X\times \Bigg(\dfrac{1}{2},1\Bigg)\rightarrow P_{cl}(X)$ be an admissible perturbation of T in the sense of Takahashi, 
  \begin{align*}
    T_W(x)=
      \begin{cases}
       \Bigg[\lambda x+(1-\lambda),~~~~~~~~\lambda x+(1-\lambda)\dfrac{1}{\sqrt{x}}\Bigg], &\text{ for } x\in \Bigg[\dfrac{1}{4},1\Bigg),\\
       ~\\
     \Bigg[\lambda x+(1-\lambda),~~~~~~~~\lambda x+(1-\lambda)\sqrt{x}\Bigg], &\text{ for } x\in [1,4].
       \end{cases}  
\end{align*}
For all $x,y\in X,$ there exists $\alpha, \beta, \gamma>0$, with $\alpha+\beta+\gamma<0$, such that $T_W$ satisfies the contraction condition of \'Ciri\'c type, as follows:
\[H(T_W(x), T_W(y))\leq \alpha d(x,y)+\beta D(x,T_W(y))+\gamma D(y, T_W(x)).\]
Then the upcoming statements hold true: 
\begin{itemize}
             \item[i)] $SFix(T)=Fix(T)=\{1\}$;
             \item[ii)] There exists $l\in \Bigg(\dfrac{1}{2\lambda},1\Bigg)$ such that 
             \[H(T(x),\{1\})\leq l H(T_W(x),\{1\}), \text{ for all } x\in X.\]
             Therefore, $(T^n(x))\overset{H}\rightarrow \{1\}$ as $n\rightarrow \infty$.
             \item[iii)] There exists $L\geq 1-\lambda $ such that \[D(x, T_W(x))\leq L D(x, T(x)), \text{ for all } x\in X.\] Then the retraction-displacement condition for T is \[|x-1|\leq \dfrac{(1+\gamma)L}{(1-\alpha-\beta)\xi}|x-1|,   \text{ for all } x\in X \text{ which is trivial for some } \xi \in (0,1).\]
\end{itemize}

\section*{Acknowledgements}
\end{example}
The author expresses her gratitude to Professor Adrian Petru\c{s}el for his support and guidance.\\

\end{document}